\documentclass{article}

\usepackage{amsfonts}
\def\cocoa
 {\mbox{\rm C\kern-.13em o\kern-.07em C\kern-.13em o\kern-.15em A}}

\newcommand\ie{{\it i.e.}}
\newcommand\eg{{\it e.g.}}

\def\round#1{\hbox{round}(#1)}
\def\openinterval#1#2{I(#1,#2)}
\def\RR{{\mathbb R}}
\def\NN{{\mathbb N}}

\title
{
Quadratic Interval Refinement for Real Roots
}

\author
{
John Abbott
}

\begin{document}
\maketitle

\begin{center}
\textbf{\Large Originally presented as a ``Poster'' at ISSAC 2006}
\end{center}

\abstract
{ We present a new algorithm for refining a real interval containing a
  single real root: the new method combines characteristics of the
  classical Bisection algorithm and Newton's Iteration.  Our method
  exhibits quadratic convergence when refining isolating intervals of simple
  roots of polynomials (and other well-behaved functions).  We assume the use
  of arbitrary precision rational arithmetic.  Unlike Newton's Iteration our
  method does not need to evaluate the derivative.
}

\section{Introduction}

Typically the process of approximating the real roots of a univariate
polynomial comprises two phases: \textit{root isolation} where the
distinct roots are separated into disjoint intervals, followed by {\it
root refinement} where the approximations of the roots are improved
until they are within specified limits.  This paper is concerned with
refinement, and assumes that isolation has already taken place (see for
instance~\cite{RouillierZimmerman}).

There are many methods for root refinement.  Newton's Iteration is a
well-known example, and dates back over three hundred years.  Bisection
is another refinement method, slower than Newton's Iteration but more
robust.  This article presents a variant of Regula Falsi which combines
the robustness of Bisection with the rapid convergence of Newton's
Iteration.  The algorithm is implemented as an integral part of {\cocoa}
(from version~$4.3$ onwards).  We call the new algorithm \textbf{
Quadratic Interval Refinement}, or simply QIR.

Prerequisites for using QIR are knowledge of an initial interval $I$ for
which the function has opposite signs at the two end points (both must
be rational), and a procedure which evaluates $f$ at a
given rational point with arbitrary precision.  Naturally, the function
$f$ must be continuous on $I$.  It is best if the interval $I$ contains
a single simple root of $f$; furthermore, if $f$ has a valid Taylor
expansion about that single simple root then convergence of QIR will
ultimately be quadratic.  If $I$ contains several roots then QIR will
eventually discard all but one of them as refinement proceeds.

\subsection{Notation}

For a real number $x$ we write $\round{x}$ to mean the integer closest
to $x$; either candidate may be chosen in the case of a tie.
We use the notation $\openinterval{a}{b}$ to denote the open interval
$\{ x \in \RR : a < x < b \}$.

We abuse terminology harmlessly by taking \textbf{isolating interval} to
mean an interval $\openinterval{a}{b}$ for which $f(a)$ and $f(b)$ have
opposite sign~---~we do not require that the interval contain just a
single root of $f$.

By \textbf{Newton's Iteration} we refer to the root approximation method
starting from some given $x_0$ and generating successive iterates using
the formula $x_n = x_{n-1} - f(x_{n-1})/f'(x_{n-1})$.  This method is
sometimes called Newton-Raphson Iteration.

By \textbf{Bisection} we refer to the interval refinement method where
the interval $\openinterval{a}{b}$ is refined to either
$\openinterval{a}{(a+b)/2}$ or $\openinterval{(a+b)/2}{b}$, the choice
depending on the sign of $f((a+b)/2)$.  This method is also called
Binary Chop.

By \textbf{Regula Falsi} we refer to the interval refinement method where
linear interpolation is used to predict the position of the root,
and then \textit{one} end of the interval is moved to this predicted
value so that the new interval still contains a sign change.

\section{Motivation}

Why invent a new method when Newton's Iteration works so well?  QIR
offers a simple unified approach which can refine any isolating
interval with a rate of convergence and a computational cost
comparable to that of Newton's Iteration.  Thus QIR has several
advantages over Newton's Iteration:
\begin{description}
\item{(a)} the criteria guaranteeing convergence are straightforward and
  easy to verify;
\item {(b)} there is no need to evaluate the derivative of the function whose
  zeroes we are approximating;
\item {(c)} the successive approximants have simple denominators;
\item {(d)} roots are approximated by intervals containing them, so
  the accuracy of an approximation is quite explicit.
\end{description}

A significant difficulty with Newton's Iteration is obtaining a suitable
starting value given an isolating interval: this is the essence of point
(a).  It may be necessary to use some other interval refinement method
prior to using Newton's Iteration.  Point (b) poses no problem when
approximating roots of explicit polynomials since the derivative can
readily be evaluated, but in other cases it may not be easy to obtain
values of the derivative.  Point~(c) is relevant only when using
arbitrary precision rational arithmetic.  Point~(d) is important when a
guarantee of the accuracy obtained is desired.  With some cunning one
can alleviate weaknesses~(c) and~(d) of Newton's Iteration:
\eg~replacing approximants by similar values having a simpler
denominator of about the right size, and using the derivative to
estimate the width of an interval containing the root.

\begin{description}
\item [Example] QIR can never fail, whereas Newton's Iteration can.  The
polynomial $f(x) = x^3-x+0.7$ has a single real root, but Newton's
Iteration fails to converge when given a starting value $x_0$ in the
range $-0.1 < x_0 < 0.1$; so bad starting points are not so rare.
\end{description}

Note that one should avoid naive use of Newton's Iteration with exact
rational arithmetic because the sizes of the numerators and
denominators of successive approximants can increase rapidly: typically
the $k$-th approximant will have numerator (or denominator) containing
$O(d^k)$ digits where $d$ is the degree of the polynomial whose
root is being sought.

\begin{description}
\item [Example] The polynomial $x^5-2$ contains a real root in the
  interval $I(1,2)$.  To estimate the root with an error less than
  $2^{-32}$ starting from the initial approximation $x_0 = 1$ requires
  five Newton iterations, and produces an approximant having a
  denominator with over $500$ digits.  Alternatively starting from
  $x_0 = 2$ requires seven iterations, and the final approximant has a
  denominator with over $17000$ digits.  In contrast, QIR obtains an
  interval of width $2^{-32}$ after six iterations, and no value in the
  computation has more than $50$ digits.
\end{description}

\section{Description of the Method}

QIR works by repeatedly narrowing the isolating interval until the width is
smaller than the prescribed limit.  Each individual narrowing step receives
an \textbf{initial interval} and a \textbf{refinement factor} $N \in \NN$
by which it must reduce the interval width; it uses \textit{discretized}
linear interpolation to guess a narrower interval containing the root.  If
the guess was good, the interval is updated; if not, the interval is
unchanged, and the narrowing step returns an indication of failure.  The
refinement factor is increased after every successful narrowing step.

A narrowing step conceptually divides its initial interval into $N$
consecutive equal width sub-intervals.  It then uses linear interpolation
to locate the sub-interval in which it guesses the root to lie.  The guess
is tested by evaluation at the end points of the chosen sub-interval.  If
the guess was lucky then the initial interval is replaced by the chosen
sub-interval; otherwise an indication of failure is returned.  The case
$N=4$ is handled specially: the initial interval is \textit{always}
replaced by the correct sub-interval, \textit{but} failure is indicated if
linear interpolation led to a bad guess.

The refinement factor is varied according to the following rule:
\begin{itemize}
\item after a successful narrowing step, $N \leftarrow N^2$
\item after a failed narrowing step, if $N > 4$ then $N \leftarrow \sqrt{N}$
\end{itemize}
This strategy is inspired by the knowledge that linear interpolation
produces approximations whose errors decrease roughly quadratically when
sufficiently close to a root of a well-behaved function.

Obviously, if one of the evaluation points happens to be the exact root
$\xi$ then this value is returned as an exact root.  This is possible
only if $\xi$ is rational and with ``suitable'' denominator.

\subsection{An Illustrative Example}

We present here an example where QIR appears to achieve quadratic
convergence immediately (but only due to ``good fortune''), then later
it is forced to reduce the refinement factor back to $4$ for a few
iterations before true quadratic convergence can begin.  Note how few
evaluations are needed to attain quite high accuracy.

\begin{description}
\item [Example] Given the polynomial $f=10^{200} \cdot x^2-1$ the root
  isolator in {\cocoa} produces the interval $\openinterval{0}{2}$
  for the positive root.  Starting from this interval QIR initially
  enjoys seven \textbf{success}es, increasing the refinement factor to
  $4^{128}$, before encountering a number of \textbf{failure}s which
  force the refinement factor back down to $4$.  Then once the
  interval width has become small enough (less than $10^{-100}$, after
  $24$ iterations in \verb|RefineInterval|) true quadratic convergence
  begins.  For instance, an interval of width less than $10^{-1000}$
  is obtained after a total of $34$ iterations, \ie~at most $68$
  evaluations of $f$; a further four iterations give an interval of
  width less than $10^{-10000}$.
\end{description}

\section{The Algorithm}

Here we present explicitly the algorithm using a pseudo-language.  The
actual source code is contained in the {\cocoa} package
\texttt{RealRoots.cpkg} and is publicly available as part of the standard
distribution of {\cocoa} (from version~$4.3$ onwards)~---~see the web
site~\cite{cocoaweb}.

\subsection{Main routine: \texttt{RefineInterval}}

\begin{description}
\item [INPUT:] $f$ the function whose zero is being sought,\\
  $\openinterval{x_{lo}}{x_{hi}}$ an open interval with rational
  end points $x_{lo} < x_{hi}$, and for which
  $f(x_{lo}) \cdot f(x_{hi})<0$\\
  $\varepsilon > 0$ an upper bound for width of the refined interval.

\item [OUTPUT:] $\openinterval{\xi_{lo}}{\xi_{hi}}$ an open
sub-interval containing a zero of $f$ and having width at most
$\varepsilon$~---~the end points are both rational.  Exceptionally the
output may be the exact value of the root $\xi$.

\item {(1)} Initialize refinement factor $N \leftarrow 4$, and
  interval $I \leftarrow \openinterval{x_{lo}}{x_{hi}}$.

\item {(2)}  While $\hbox{\it width}(I) > \varepsilon$ do steps (2.1)--(2.4).

\begin{description}
\item {(2.1)} Apply \texttt{RefineIntervalByFactor} to $f$, $I$ and
  $N$~---~this will usually modify the interval $I$.
\item{(2.2)}  If an \textbf{exact root} $\xi$ has been found, simply return $\xi$.
\item{(2.3)}  If \textbf{success} was reported then replace $N \leftarrow N^2$.
\item{(2.4)}  Otherwise \textbf{failure} was reported, so if $N > 4$ replace $N \leftarrow \sqrt{N}$.
\end{description}

\item {(3)} Return the interval $I$.
\end{description}

\subsection{Auxiliary procedure:  \texttt{RefineIntervalByFactor} when $N > 4$}

\begin{description}
\item [INPUTS:] $f$ the function whose zero we are seeking,\\
  $I = \openinterval{x_{lo}}{x_{hi}}$ an open interval with
  rational end points $x_{lo} < x_{hi}$, and for which
  $f(x_{lo}) \cdot f(x_{hi})<0$\\
  $N$ the refinement factor, and we assume $N>4$~---~see \S\ref{case-N-equals-4} for the case $N=4$.

\item [OUTPUT:] \textbf{success} or \textbf{failure}; if \textbf{success} then as a
 side effect the interval $I$ is replaced by a sub-interval whose
 width is $1/N$ times that of $I$ and which contains a zero of $f$.
 Exceptionally the output may be \textbf{exact root} $\xi$.

\item {(1)} Conceptually divide the interval $I$ into $N$ consecutive
  sub-intervals of equal width.  Compute the width of each sub-interval:
  $w = (x_{hi} - x_{lo})/N$.  Define $x_0 = x_{lo}$, and $x_k =
  x_{k-1}+w$ for $k=1,2,\ldots,n$; so these $x_i$ are the end points of
  the sub-intervals.

\item {(2)} Using linear interpolation predict approximate position of
  the zero: \ie~determine the closest end point $\hat x = x_{\kappa}$
  where $\kappa = \round{N \cdot f(x_{lo}) / (f(x_{lo}) - f(x_{hi}))}$

\item {(3)} If $f(\hat x) = 0$ then return \textbf{exact root} $\hat x$.

\item {(4)} Check whether our prediction was good; \ie~the
  function $f$ changes sign in one of the sub-intervals having $\hat
  x$ as an end point.

\item {(5)} If $f(\hat x)$ has the same sign as $f(x_{lo})$, check
  the interval to the right:
\begin{description}
\item {(5.1)}  If $f(\hat x + w) = 0$ then return  \textbf{exact root} $\hat x +w$.
\item {(5.2)}  If $f(\hat x +w)$ has the same sign as
  $f(\hat x)$ then return \textbf{failure}.
\item {(5.3)} Otherwise replace $I \leftarrow \openinterval{\hat x}{\hat x + w}$ and return \textbf{success}.
\end{description}

\item {(6)} Otherwise $f(\hat x)$ has the same sign as $f(x_{hi})$, so check the interval to the left:
\begin{description}
\item  If $f(\hat x - w) = 0$ then return \textbf{exact root} $\hat x - w$.
\item  If $f(\hat x - w)$ has the same sign as $f(\hat x)$ then return \textbf{failure},
\item  Otherwise replace $I \leftarrow \openinterval{\hat x - w}{\hat x}$  and return \textbf{success}.
\end{description}
\end{description}

\subsection{Auxiliary procedure:  \texttt{RefineIntervalByFactor} when $N=4$}
\label{case-N-equals-4}

We give a verbal description rather than pseudo-code for
\texttt{RefineIntervalByFactor} with $N=4$.  Use linear interpolation to
predict which of the $5$ sub-interval end points is closest to the root.
Refine the input interval using Bisection twice.  If one of the end points
of the refined interval agrees with our prediction, return
\textbf{success}, otherwise return \textbf{failure}.  Note that the initial
interval is \textit{always} narrowed by a factor of $4$.


\subsection{Termination of the Algorithm}

We show that the algorithm always terminates.  In \texttt{RefineInterval}
if $N=4$ then the interval is always narrowed; if $N>4$ then either the
interval is narrowed (by more than a factor of $4$), or $N$ is reduced.
Since \texttt{RefineInterval} terminates when the interval has width less
than $\varepsilon$, and since each narrowing reduces the width by at
least a factor of $4$, there can be only a finite number of narrowings.
If no narrowing occurs (\ie~$N>4$ and \textbf{failure} was reported) then
$N$ is reduced, but $N$ is never smaller than $4$, so there can be only
finitely many iterations in which the interval is not narrowed.

\subsection{Quadratic Convergence}\label{quadratic}

The mathematical justification underlying the convergence rate of QIR is
very similar to that underlying Newton's Iteration.  We assume that the
isolating interval $I$ contains a single simple root $\xi$ and that the
function $f$ admits a Taylor expansion centred on $\xi$ valid in an open
neighbourhood: \ie~for sufficiently small $\delta$ we have
$$
f(\xi+ \delta) =  C_1 \cdot \delta + C_2 \cdot \delta^k + O(\delta^{k+1})
$$
for some exponent $k \ge 2$ and suitable constants $C_1$ and $C_2$; 
moreover $C_1 \neq 0$ since $\xi$ is a simple root.

Let $x_{lo}$ and $x_{hi}$ be the two end points of the interval $I$ and
$\varepsilon$ be its width.  When $\varepsilon$ is small enough we have:
\begin{eqnarray*}
f(x_{lo}) &=& C_1 \cdot (x_{lo}-\xi) + C_2 \cdot (x_{lo}-\xi)^k + O(x_{lo}-\xi)^{k+1}\\
f(x_{hi}) &=& C_1 \cdot (x_{hi}-\xi) + C_2 \cdot (x_{hi}-\xi)^k + O(x_{hi}-\xi)^{k+1}
\end{eqnarray*}
Estimating $\xi$ by
linear interpolation gives $\hat\xi = x_{lo} + \lambda \cdot (x_{hi} -
x_{lo})$ where $\lambda = f(x_{lo})/(f(x_{lo})-f(x_{hi}))$.
Observe that $f(x_{lo})-f(x_{hi}) = C_1 \varepsilon + C_2 \eta +O(\varepsilon^{k+1})$ where $|\eta| = |(x_{lo}-\xi)^k - (x_{hi}-\xi)^k| \le \varepsilon^k$.
Substituting and simplifying gives $\hat\xi = x_{lo} + (\xi-x_{lo}) - \frac{C_2}{C_1} (x_{lo}-\xi)( \frac{\eta}{\varepsilon} - (x_{lo}-\xi)^{k-1}) + O(\varepsilon^{k+1})$.
Regarding cubic and higher order terms in $\varepsilon$ as negligible, this
 tells us that $|\hat\xi - \xi| < \varepsilon^2 \cdot 2C_2/C_1$.


\medskip

Now we show that convergence of QIR is eventually quadratic.  We say that
the isolating interval $I$ of width $\varepsilon$ is \textbf{narrow enough}
for $f$ whenever $\varepsilon \cdot C_2/C_1 < 1/8$.  Observe that, starting
from any isolating interval, this condition will be satisfied after only
finitely many iterations.  Henceforth we shall assume $I$ is narrow enough.

When applying the algorithm \texttt{RefineIntervalByFactor} to the
parameters $(f, I, N)$ we shall say that $N$ is \textbf{small enough for}
$I$ if $N \cdot \varepsilon\cdot C_2/C_1 < 1/2$.  Note that $N=4$ is always
small enough because we have assumed that $I$ is narrow enough.

Suppose that algorithm \texttt{RefineIntervalByFactor} is called on $(f,
I, N)$ where $N$ is small enough for $I$.  The condition
$\varepsilon \cdot C_2/C_1 < 1/8$ implies that linear interpolation produces an estimate
$\hat\xi$ whose distance from $\xi$ does not exceed $\varepsilon^2 \cdot
2C_2/C_1 < \varepsilon/N$.  Hence step~(2) of
\texttt{RefineIntervalByFactor} makes a good prediction, and so
\textbf{success} will be returned.  Denote the narrowed interval by $I'$.
Now, the next iteration of \texttt{RefineInterval} will call
\texttt{RefineIntervalByFactor} on the parameters $(f, I', N^2)$.
Since $I'$ is narrower than $I$ it is surely narrow enough for $f$.  Moreover,
one may easily verify that $N^2$ is small enough for $I'$.  Hence, by
induction, all subsequent calls to \texttt{RefineIntervalByFactor} will
return \textbf{success}, and convergence is therefore quadratic.

In contrast, while still assuming that $I$ is narrow enough for $f$, if
\texttt{RefineIntervalByFactor} is called when $N$ is not small enough for
$I$ then \textbf{failure} may occur.  Since $N$ is not small enough, we must
have $N>4$.  So when \textbf{failure} occurs the refinement factor $N$ is
reduced.  Hence, there can be only finitely many \textbf{failure}s before $N$
becomes small enough for $I$, and then guaranteed quadratic convergence ensues.

\subsection{Weaknesses}

A weakness of QIR as described here is the cost of evaluating
$f$ exactly at various rational points.  For instance, if $f$ is a high
degree polynomial and the evaluation point $x$ is a rational with large
denominator then $f(x)$ is likely to have an enormous denominator:
more or less the $n$-th power of the denominator of $x$ where $n$ is the
degree of $f$.

An interesting solution would be to use high precision floating point
arithmetic.  If variable precision floating point numbers are available
then it would be sufficient to evaluate $f$ approximately to obtain values
with at least $m+2$ bits of precision where $m=\log_2(N)$ the logarithm of
the refinement factor.  If the relative errors in $f(x_{lo})$ and
$f(x_{hi})$ are at most $2^{-m-2}$ then step~(2) in
\texttt{RefineIntervalByFactor} will compute the index $\kappa$ differing
by at most $1$ from the value one would have obtained using exact rational
arithmetic.

\medskip
As presented here QIR may calculate an interval much narrower than the
specified limit. This is a (minor) defect because it means the last two
evaluations of $f$ will be at rational numbers with needlessly large
denominators (with up to twice as many digits as truly required), and
thus costing more than is strictly necessary.  In practice one can
simply reduce the value of the refinement factor $N$ in the final
iteration of {\verb|RefineInterval|} to avoid the overshoot.

\section{Experimental Results}

We present here a small selection of experimental results.  The computer
used was a Macintosh G5 (running MacOS~X~10.4) with 2GHz processors, and
3.5G bytes of RAM.  The software used was {\cocoa~4.6} which relies upon
the excellent library GMP (version~$4.2$, see~\cite{GMPweb}) for
arbitrary precision arithmetic.  All timings are measured in seconds.

The real root isolation code in {\cocoa} is a \textit{simplistic}
implementation of the method described in \cite{RouillierZimmerman},
which is based on Descartes's Rule of Signs.  Their implementation
appears to be faster than ours in {\cocoa}.

We describe briefly the polynomials used for the timing tests.  The
first two polynomials are adapted from the FRISCO test suite
(see~\cite{friscoweb1} and~\cite{friscoweb2}) since {\cocoa} cannot
represent complex numbers.  Both polynomials are designed to be
challenging for a root isolator, and indeed the real root isolator implemented
in {\cocoa} fails with \textit{``Recursion too deep''}.  However
$\openinterval{0}{1}$ is a valid isolating interval for both
polynomials, and it was used as input to QIR.  The first polynomial is
$f_1 = ((c^2x^2-3)^4+c^4x^{18})(c^2x^2-3)$ for $c=10^{100}$; there are
four complex roots quite close to the real positive root, and this
delays the onset of quadratic convergence.  The second polynomial is
$f_2 = x^{50}+(10^{50}x-1)^3$ and has two complex roots very close
to the real positive root again delaying the onset of quadratic
convergence: note how quickly $10000$ digits are obtained compared to
the time needed for $1000$ digits.

The third and fourth polynomials have a similar structure; all their
roots are real by construction.  $f_3$ has
degree $32$ and is the minimal polynomial of
$\sqrt{176}+\sqrt{190}+\sqrt{195}+\sqrt{398}+\sqrt{1482}$, while $f_4$
has degree $128$ and is the minimal polynomial of
$\sqrt{99}+\sqrt{627}+\sqrt{661}+\sqrt{778}+\sqrt{929}+\sqrt{1366}+\sqrt{1992}$.
These polynomials are interesting because they both have a root
which is ``almost integer'': $f_3$ has a root close to $45$, and $f_4$
has a root close to $10$.

In the table below, the column headed ``\#roots'' indicates how many
intervals were refined, the column headed ``isolate'' is the time
taken by the root isolator, and the other columns give the times
taken to refine all the intervals found by the isolator to a width not
greater than $\varepsilon$.

\smallskip\noindent
\begin{center}
\begin{tabular}{|l|r||r|r|r|r|}
\hline
Poly  &\#roots&isolate&$\varepsilon=10^{-100}$&$\varepsilon=10^{-1000}$&$\varepsilon=10^{-10000}$
$\vphantom{\int_0^0}$\\
\hline
$f_1$ &    1 &    fail&     0.0  &          2.4 &          2.9\\
$f_2$ &    1 &    fail&     0.4  &         43~~ &         47~~\\
$f_3$ &   32 &    2.8 &     0.6  &          2.4 &        100~~\\
$f_4$ &  128 &   74~~ &    10~~  &        125~~ &       6600~~\\
\hline
\end{tabular}
\end{center}

\bigskip
The GMP library offers functions for computing the integer part of roots of
big integers (\textit{viz.}~\verb|mpz_sqrt| and \verb|mpz_root|); these
functions use Newton's Iteration specialized to the computation of $k$-th
roots.  We compared QIR with these functions in the following way.  To
compute $d$ digits of $\sqrt{5}$ using \verb|mpz_sqrt| we passed as
argument the value $5 \times 10^{2d}$.  Similarly to compute $d$ digits of
$\root{3}\of{3}$ and $\root{5}\of{2}$ we called \verb|mpz_root| with
arguments $3 \times 10^{3d}$ and $2 \times 10^{5d}$ respectively.  The
timings for QIR are for refining the root interval to a width not exceeding
$10^{-d}$~---~the initial isolating interval was $[0,4]$ in both cases.
The square root function \verb|mpz_sqrt| is especially efficient: \eg~it
computed $10^7$ digits of $\sqrt{5}$ in $8.2$~seconds while QIR needed
$35$~seconds to achieve the same accuracy.  In contrast, we see that QIR is
quite competitive with \verb|mpz_root|; here are the results:

\smallskip\noindent
\begin{center}
\begin{tabular}{|c||c|c||c|c|}
\hline
                   & QIR      &  QIR     & GMP      & GMP \\
Value approximated & $d=10^5$ & $d=10^6$ & $d=10^5$ & $d=10^6$\\
\hline
$\root{3}\of{3}$   &   0.4    &  5.1     &   0.8    &  17 \\
$\root{5}\of{2}$   &   1.1    & 12~~~~   &   1.3    &  26 \\
\hline
\end{tabular}
\end{center}

\section{Conclusions}

QIR is a simple algorithm combining the speed of Newton's Iteration with
the robustness of Bisection.  It is easy to implement, and the
implementation supplied with {\cocoa} demonstrates its genuine
competitiveness with Newton's Iteration (except perhaps for computing
square roots).


\end{document}